\documentclass[11pt]{article}
\usepackage{latexsym}
\usepackage{graphicx}
\usepackage{amsmath}
\usepackage{amssymb}
\usepackage{color}
\usepackage{amscd}
\usepackage{epsfig}
\usepackage{cite}
\newtheorem{problem}{Problem}[section]

\newtheorem{lemma}[problem]{Lemma}
\newtheorem{theorem}[problem]{Theorem}
\newtheorem{prop}[problem]{Proposition}
\newtheorem{corollary}[problem]{Corollary}

\title{The $L$-functions of Witt coverings}
\author{Chunlei Liu and Dasheng Wei
\\
\small {Beijing Normal University, Beijing 100875. {\it Email:} clliu@bnu.edu.cn}\\
\small{University of Science and Technology of China, Hefei
230026. {\it Email:} dshwei@ustc.edu}}
\begin{document}
\maketitle
{\bf Abstract. }Results on $L$-functions of Artin-Schreier
coverings by Dwork, Bombieri and Adolphson-Sperber are generalized
to $L$-functions of Witt coverings.


\section{Introduction}
We shall state our main results after recalling the notion of
$L$-functions of Witt coverings.

Let $\mathbb{F}_q$ be the finite field of characteristic $p$ with
$q$ elements, and $W_m$ the ring scheme of Witt vectors of length
$m$ over $\mathbb{F}_q$. Let $f\in
W_m(\mathbb{F}_q[x_1^{\pm1},\cdots,x_n^{\pm1}])$ with its first
coordinate non-constant. Let $T^n$ be the $n$-dimensional toruse
over $\mathbb{F}_q$, and $F$ the Frobenius morphism of $W_m$. The
fibre product over $W_m$ of $W_m\stackrel{F-1}{\rightarrow}W_m$
and $T^n\stackrel{f}{\rightarrow}W_m$ is a
$W_m(\mathbb{F}_p)$-covering of $T^n$, with group action
$g(y,x)=(y+g,x)$. The Frobenius element of the Galois group
$W_m(\mathbb{F}_p)$ at a closed $x$ of $X$ with degree $k$ is
$\text{Tr}_{W_m(\mathbb{F}_{q^k})/W_m(\mathbb{F}_p)}(f(x))$. So
the Artin $L$-function of $T^n$ determined by that
$W_m(\mathbb{F}_p)$-covering and a fixed character
$\psi:W_m(\mathbb{F}_p)\rightarrow\overline{\mathbb{Q}}^{\times}$
of exact order $p^m$ is
$$L_f(t)=\prod\limits_{x\in |T^n|}
(1-\psi({\text
Tr}_{W_m(\mathbb{F}_{q^k})/W_m(\mathbb{F}_p)}(f(x)))t)^{(-1)^n},$$
where $|T^n|$ is the set of closed points of $T^n$. By a well
known theorem of Deligne [De],
$$L_f(t)=\frac{\prod\limits_{\alpha}(1-\alpha t)}{\prod\limits_{\beta}(1-\beta t)},$$
where $\alpha$ and $\beta$ are algebraic integers such that
$q^n\alpha^{-1}$ and $q^n\beta^{-1}$ are also algebraic integers.
It implies, as observed by Bombieri [Bo2],  $\text{ord}_q(\alpha),
\text{ord}_q(\beta)\leq n$, where $\text{ord}_q$ is the $q$-order
function of $\overline{\mathbb{Q}}_p$ such that
$\text{ord}_q(q)=1$. ($\overline{\mathbb{Q}}_p$ is the algebraic
closure of $\mathbb{Q}_p$, the field of $p$-adic numbers.)

By logarithmic differentiation, we get
$$L_f(t)=\exp(\sum\limits_{k=1}^\infty
S_k(f)\frac{t^k}{k}),$$where
$$S_k(f)=(-1)^{n-1}\sum\limits_{x\in (\mathbb{F}_{q^k}^{\times})^n} \psi({\text
Tr}_{W_m(\mathbb{F}_{q^k})/W_m(\mathbb{F}_p)}(f(x)))$$ are
exponential sums associated to characters of $p$-power order. To
have a look at these exponential sums, we denote by
$\lambda_i:A^1\rightarrow W_m$, $i=0,\cdots,m-1$, the embedding
which maps $A^1$ onto the $i$-th axis of $W_m$, and write $$f=\sum
\limits_{i=0}^{m-1}\sum\limits_{u\in I_i}\lambda_i(a_{iu}x^u),$$
where $I_i\subset\mathbb{Z}^n$ and
$a_{iu}\in\mathbb{F}_q^{\times}$ are uniquely determined. That
decomposition can be obtained by solving the congruences
$$\begin{array}{c}
   f\equiv \lambda_0(\sum\limits_ua_{0u}x^u)(\text{mod }V)\\
  f-\sum\limits_u\lambda_0(a_{0u}x^u)
\equiv \lambda_1(\sum\limits_ua_{1u}x^u)
(\text{mod }V^2) \\
  \vdots \\
  f-\sum\limits_{i=0}^{m-2}\sum\limits_u\lambda_i(a_{0u}x^u)
\equiv \lambda_{m-1}(\sum\limits_ua_{(m-1)u}x^u)
(\text{mod }V^m)  \\
\end{array}$$ successively, where $V$ is the shift operator on
$W_m$.

Let $\overline{\mathbb{F}}_q$ be the algebraic closure of
$\mathbb{F}_q$, and $\omega$ the Teichm\"{u}ller lifting from
$\overline{\mathbb{F}}_q$ to $\overline{\mathbb{Q}}_p$. We define
$\omega(f)=\sum \limits_{i=0}^{m-1}p^i\sum\limits_{u\in
I_i}\omega(a_{iu})x^u$. Let $\mathbb{Z}_p$ be the ring of $p$-adic
integers, and $\mu_l$ ($l\geq1$) be the set of $l$-th roots of
unity in $\overline{\mathbb{Q}}_p$. Identifying
$W_m(\mathbb{F}_{q^k})$ with $\mathbb{Z}_p[\mu_{q^k-1}]/(p^m)$
under the isomorphism
$$(a_0,\cdots,a_{m-1})\mapsto
\sum\limits_{j=0}^{m-1}\omega(a_i^{p^{-i}})p^i\ (\text{mod
}p^m),$$ one finds, for $x\in (\mathbb{F}_{q^k}^{\times})^n$, that
$$\psi({\text
Tr}_{W_m(\mathbb{F}_{q^k})/W_m(\mathbb{F}_p)}(f(x)))=\psi({\text
Tr}_{\mathbb{Q}_p[\mu_{q^k-1}]/\mathbb{Q}_p}(\sum\limits_{i=0}^{m-1}\sum\limits_{u}
p^i\omega(a_{iu}^{p^{-i}}x^{p^{-i}u})))$$$$= \psi({\text
Tr}_{\mathbb{Q}_p[\mu_{q^k-1}]/\mathbb{Q}_p}(\sum\limits_{i=0}^{m-1}\sum\limits_{u}
p^i\omega(a_{iu}x^u)))= \psi({\text
Tr}_{\mathbb{Q}_p[\mu_{q^k-1}]/\mathbb{Q}_p}(\omega(f)(\omega(x)))).$$
Therefore, we have
\begin{lemma}For $k=1,2,\cdots$, we have
$$S_k(f) =\sum\limits_{x\in\mu_{q^k-1}^n}
\psi({\text
Tr}_{\mathbb{Q}_p[\mu_{q^k-1}]/\mathbb{Q}_p}(\omega(f)(x))).$$
\end{lemma}

We define the Newton polyhedron $\Delta_{\infty}(f)$ of $f$ at
infinity to be the convex hull in $\mathbb{Q}^n$ of
$\{p^{m-i-1}u:0\leq i\leq m-1,u\in I_i\}\cup\{0\}$. Recall that,
for a convex polyhedron $\Delta$ of dimension $n$ in
$\mathbb{Q}^n$ that contains the origin, there is a
$\mathbb{R}_{\geq0}$-linear degree function $u\mapsto\deg(u)$ on
$L(\Delta)$, the set of integral points in the cone
$\bigcup\limits_{k=1}^{\infty}k\Delta$, such that $\deg(u)=1$ when
$u$ lies on a face of $\Delta$ that does not contain the origin.
That degree function may take on non-integral values. But there is
a positive integer $D$ such that $\deg L(\Delta)\subset
D^{-1}\mathbb{Z}$. We denote the least positive integer with this
property by $D(\Delta)$. For $k=0,1,\cdots$, we denote by
$W_{\Delta}(k)$ the number of points of degree
$\frac{k}{D(\Delta)}$ in $L(\Delta)$. We define
$P_{\Delta}(t)=(1-t^{D(\Delta)})^n\sum\limits_{k=0}^{+\infty}W_{\Delta}(k)t^k$
for later use. Our first result is an upper bound for the total
degree of $L_f(t)$.
\begin{theorem}The total degree of $L_f(t)$ is bounded by
$\sum\limits_{i=0}^n(\begin{array}{c}
  n \\
  i \\
\end{array})\sum\limits_{k=0}^{D(n-i+1)}W_{\Delta}(k)$
with $D=D(\Delta)$ and $\Delta=\Delta_{\infty}(f)$.\end{theorem}

For $j=1,\cdots,n$, we write
$$\overline{_jf}^{\tau}=\sum
\limits_{i=0}^{m-1}\sum\limits_{p^{m-i-1}u\in\tau}u_ja_{iu}^{p^{m-i-1}}x^{p^{m-i-1}u},$$
where $u_j$ is the $j$-th coordinate of $u$. We call $f$
non-degenerate with respect to $\Delta_{\infty}(f)$ if
$\Delta_{\infty}(f)$ is of dimension $n$, and for every face
$\tau$ of $\Delta_{\infty}(f)$ that does not contain 0, the system
$\overline{_1f}^{\tau} =\cdots=\overline{_nf}^{\tau}$ has no
common solution in $(\overline{\mathbb{F}}_q^{\times})^n$. Our
second result is on $L$-functions from non-degenerate Witt
vectors.
\begin{theorem}Suppose that $f$ is non-degenerate with respect to
$\Delta:=\Delta_{\infty}(f)$. Then the $L$-function $L_f(t)$ is a
polynomial, and its Newton polygon with respect to $\text{ord}_q$
lies above the Hodge polygon of $P_{\Delta}(t)$ of degree
$D(\Delta)$ with the same endpoints. In particular, $L_f(t)$ is of
degree $n!\text{Vol}(\Delta)$.
\end{theorem}
Recall that the Newton polygon of $\prod(1-\alpha
t)\in\overline{\mathbb{Q}}_p[[t]]$ with respect to $\text{ord}_q$
is the polygon with vertices at points
$$(\sum\limits_{\text{ord}_q(\alpha)\leq y}1,
\sum\limits_{\text{ord}_q(\alpha)\leq y}\text{ord}_q(\alpha)), \
y\in\mathbb{Q}.$$ And the Hodge polygon of
$\sum\limits_{k=0}^{+\infty}a_kt^k$ of degree $D$ is the polygon
with vertices at the points $(0,0)$ and
$$(\sum\limits_{i=0}^ka_i,\frac{1}{D}\sum\limits_{i=0}^kia_i),
\ k=0,1,\cdots.$$

Theorem 1.3 was proved by Dwork [Dw] when $m=1$, and
$f(x_1,\cdots,x_n)=x_nh(x_1,\cdots,x_{n-1})$ for some polynomial
$h$ with coefficients in $\mathbb{F}_q$. In that case, the
$L$-function $L_f(t)$, by the orthogonality of characters, is
related to the zeta function of the hypersurface defined by $h=0$
in the $(n-1)$-dimensional affine space defined over
$\mathbb{F}_q$. It was completely proved by Adolphson-Sperber
[AS2] in the case $m=1$. In the case $n=1$, the degree of $L_f(t)$
was determined by Kumar-Helleseth-Calderbank [KHC] with
applications to coding theory, and by W.-C. W. Li [Li], who read
the $p=2$ version of [KHC].

Our proof of the main results is based on the $p$-adic method set
up by Dwork [Dw, Dw2] and developed by Bombieri [Bo, Bo2], Monsky
[Mo], Adolphson-Sperber [AS, AS2], Wan [Wn], and others. The
innovation lies in the use of the Artin-Hasse exponential series
to produce roots of unity of $p$-power order.

One can infer the following theorem from Theorem 1.3.
\begin{theorem}If $f$
is non-degenerate with respect to $\Delta_{\infty}(f)$, and the
origin lies in the interior of $\Delta_{\infty}(f)$, then the
reciprocal roots of $L_f(t)$ are of absolue value $q^{n/2}$.
\end{theorem}

{\bf Acknowledgement.} This work is completed when the authors are
visiting the Morningside Center of Mathematics in Beijing (MCM).
The authors thank MCM for its hospitality, and Fei Xu for inviting
them to visit MCM. The authors also want to express their
gratitude to Daqing Wan, who teaches them $p$-adic analysis,
encourages them to study exponential sums, and patiently guides
their research. Chunlei Liu would like to thank Yingbo Zhang for
support, and Lei Fu for discussions. The research of Chunlei Liu
is supported by NSFC Grant No. 10371132, by Project 985 of Beijing
Normal University, and by the Foundation of Henan Province for
Outstanding Youth.
\section{The Artin-Hasse exponential series} Let
$$E(t)=\exp(\sum_{i=0}^{\infty}\frac{t^{p^i}}{p^i}) \in
\mathbb{Z}_p[[t]]$$ be the Artin-Hasse exponential series. We
shall use it to produce roots of unity of $p$-power order.
\begin{lemma}If $l$ is a positive integer, and $\pi$ is
 a root of $\sum\limits_{i=0}^{\infty}\frac{t^{p^i}}{p^i}=0$ in
$\overline{\mathbb{Q}}_p$ with order $\frac{1}{p^{l-1}(p-1)}$,
$E(\pi)$ is a primitive $p^l$-th root of unity.
\end{lemma}
{\it Proof.  } First, $\exp(p^l\frac{\pi^{p^i}}{p^i})$ exists as
$\text{ord}_p(p^l\frac{\pi^{p^i}}{p^i})\geq\frac{p}{p-1}$. So
$$E(\pi)^{p^l}=E(p^lt)|_{t=\pi}=
\prod\limits_{i=0}^{\infty}\exp(p^l\frac{\pi^{p^i}}{p^i})
=\exp(\sum\limits_{i=0}^{\infty}p^l\frac{\pi^{p^i}}{p^i})=\exp(0)=1.$$
Secondly, as $E(t)\in1+t+t^2\mathbb{Z}_p[[t]]$,
$$E(\pi)^{p^{l-1}}\equiv(1+\pi)^{p^{l-1}}\equiv1+\pi^{p^{l-1}} (\text{ mod }\pi^{p^{l-1}+1}).$$
The lemma is proved.
\begin{lemma}Let $l$ be a positive integer.
Then the Artin-Hasse exponential series induces a bijection $\pi\mapsto
E(\pi)$ from the set of roots of
$\sum\limits_{i=0}^{\infty}\frac{t^{p^i}}{p^i}=0$ in
$\overline{\mathbb{Q}}_p$ with order $\frac{1}{p^{l-1}(p-1)}$ to
the set of all primitive $p^l$-th roots of unity in
$\overline{\mathbb{Q}}_p$.
\end{lemma}
{\it Proof.} The field generated over $\mathbb{Q}_p$ by the set of
roots of $\sum\limits_{i=0}^{\infty}\frac{t^{p^i}}{p^i}=0$ in
$\overline{\mathbb{Q}}_p$ with order $\frac{1}{p^{l-1}(p-1)}$ is
precisely $\mathbb{Q}_p(\mu_{p^l})$ since it contains
$\mathbb{Q}_p(\mu_{p^l})$ by the preceding lemma, and is of degree
no greater than $p^{l-1}(p-1)$ over $\mathbb{Q}_p$ by Weierstrass'
Preparation Theorem. One sees that $E(\tau(\pi))=\tau(E(\pi))$ if
$\tau$ is an automorphism $\mathbb{Q}_p(\mu_{p^l})$ over
$\mathbb{Q}_p$. So $\pi\mapsto E(\pi)$ maps the set of roots of
$\sum\limits_{i=0}^{\infty}\frac{t^{p^i}}{p^i}=0$ in
$\overline{\mathbb{Q}}_p$ with order $\frac{1}{p^{l-1}(p-1)}$ onto
the set of all primitive $p^l$-th roots of unity in
$\overline{\mathbb{Q}}_p$. It is a bijection as
$\sum\limits_{i=0}^{\infty}\frac{t^{p^i}}{p^i}=0$ has at most
$p^{l-1}(p-1)$ roots in $\overline{\mathbb{Q}}_p$ with order
$\frac{1}{p^{l-1}(p-1)}$ by Weierstrass' Preparation Theorem.
\begin{lemma}If
$k$ is a positive integer, and $x \in \overline{\mathbb Q}_p$
satisfies $x^{p^k}=x$, then
$$E(t)^{x+x^p+\cdots+x^{p^{k-1}}}=E(tx)E(tx^p)\cdots E(tx^{p^{k-1}}).$$
\end{lemma}
{\it Proof. } As
$\sum\limits_{j=0}^{k-1}x^{p^j}=\sum\limits_{j=0}^{k-1}x^{p^{j+i}}$,
we have
$$
E(t)^{x+x^p+\cdots+x^{p^{k-1}}}
=\exp(\sum_{i=0}^{\infty}\frac{t^{p^i}}{p^i}\sum\limits_{j=0}^{k-1}x^{p^j})
=\exp(\sum_{i=0}^{\infty}\frac{t^{p^i}}{p^i}\sum\limits_{j=0}^{k-1}x^{p^{j+i}})$$
$$
=\exp(\sum\limits_{j=0}^{k-1}\sum_{i=0}^{\infty}\frac{(tx^{p^j})^{p^i}}{p^i})=E(t
x)E(t x^p)\cdots E(t x^{p^{k-1}}).$$ The lemma is proved.
\begin{corollary}If $\pi$ is
 a root of $\sum\limits_{i=0}^{\infty}\frac{t^{p^i}}{p^i}=0$ in
$\overline{\mathbb{Q}}_p$ with order $\frac{1}{p^{l-1}(p-1)}$, and
$x \in \overline{\mathbb Q}_p$ satisfies $x^{p^k}=x$, then
$$E(\pi)^{x+x^p+\cdots+x^{p^{k-1}}}=E(\pi x)E(\pi x^p)\cdots E(\pi x^{p^{k-1}}).$$
\end{corollary}

We now fix an embedding of $\overline{\mathbb{Q}}$ into
$\overline{\mathbb{Q}}_p$. Guaranteed by the above lemma, we may
choose, for each $l=1,\cdots,m$, a unique root $\pi_l$ of
$\sum\limits_{i=0}^{\infty}\frac{t^{p^i}}{p^i}=0$ in
$\overline{\mathbb{Q}}_p$ with order $\frac{1}{p^{l-1}(p-1)}$ such
that $E(\pi_l)=\psi(1)^{p^{m-l}}$. Let
$\Delta=\Delta_{\infty}(f)$, $D=D(\Delta)$, and $\pi$ a $D$-th
root of $\pi_m^{p^{m-1}}$ in $\overline{\mathbb{Q}}_p$. For
$b\geq0$, we write
$$L(b)=\{\sum_{u\in L(\Delta)}a_ux^u :
 a_u\in \mathbb{Z}_p[\mu_{q-1},\pi_m,\pi], \text{ ord}_p(a_u)\geq b\deg(u)\}.$$
The Galois group
$\text{Gal}(\mathbb{Q}_p[\mu_{q-1},\pi_m,\pi]/\mathbb{Q}_p)$ acts
on $L(b)$ coefficientwise. Define $$E_f(x)
=\prod\limits_{i=0}^{m-1}\prod\limits_{u\in
I_i}E(\pi_{m-i}\omega(a_{iu})x^u).$$
\begin{lemma}We have $E_f(x)\in L(\frac{1}{p-1})$.\end{lemma}
{\it Proof. }Suppose that $0\leq i\leq m-1$ and $u\in I_i$. We
have $p^{m-i-1}u\in\Delta$. So $\deg(p^{m-i-1}u)\leq1$, and
$$\text{ord}_p(\pi_{m-i})=\frac{1}{p^{m-i-1}(p-1)}\geq
\frac{\deg(p^{m-i-1}u)}{p^{m-i-1}(p-1)}=\frac{\deg(u)}{p-1}.$$ It
follows that $\pi_{m-i}\omega(a_{iu})x^u\in L(\frac{1}{p-1})$.
Since $E(t)\in\mathbb{Z}_p[[t]]$, we have
$E(\pi_{m-i}\omega(a_{iu})x^u)\in L(\frac{1}{p-1})$. The lemma now
follows.

Let $\sigma$ be the Frobenius element of
$\text{Gal}(\mathbb{Q}_p[\mu_{q-1},\pi_m,\pi]/\mathbb{Q}_p)$
fixing $\pi_m$ and $\pi$. The following lemma follows from
Corollary 2.4.
\begin{lemma}If
$k$ is a positive integer, and $x \in \mu_{q^k-1}^n$, then
$$\psi({\text
Tr}_{\mathbb{Q}_p[\mu_{q^k-1}]/\mathbb{Q}_p)}(\omega(f)(x)))
=\prod\limits_{i=0}^{ak-1}E_f^{\sigma^i}(x^{p^i}).$$
\end{lemma}
\begin{corollary}We have
$$S_k(f)=(-1)^{n-1}\sum\limits_{x\in\mu_{q^k-1}^n}
\prod\limits_{i=0}^{ak-1}E_f^{\sigma^i}(x^{p^i}),\ k=1,2,\cdots.$$
\end{corollary}
\section{Functions from the Artin-Hasse exponential series}
We shall study the growth of the coefficients of $\widehat{_kf}$
($k=1,\cdots,n$), which are defined by
$$d\log\widehat{E}_f(x)=
\sum\limits_{k=1}^n\widehat{_kf}\frac{dx_k}{x_k},\
\widehat{E}_f(x)=\prod\limits_{j=0}^{\infty}
E_f^{\sigma^j}(x^{p^j}).$$

\begin{lemma}We have$$\widehat{_kf}=
\sum\limits_{i=0}^{m-1}\sum\limits_{j=0}^{\infty}p^j\gamma_{i,j}\sum\limits_{u\in
I_i}u_k\omega(a_{iu}^{p^j})x^{p^ju},\ k=1,\cdots,n,$$ where
$\gamma_{i,j}=\sum\limits_{l=0}^j\frac{\pi_{m-i}^{p^l}}{p^l}$.\end{lemma}
\begin{lemma}We have
$\pi_l\equiv\pi_m^{p^{m-l}} (\text{ mod
}\pi_m^{p^{m-l}+1})$.\end{lemma}Since
$E(t)\in1+t+t^2\mathbb{Z}_p[[t]]$, we have $E(\pi_l)\equiv1+\pi_l
(\text{ mod }\pi_l^2)$. So we have
$$E(\pi_m)^{p^{m-l}}\equiv(1+\pi_m)^{p^{m-l}}\equiv1+\pi_m^{p^{m-l}}
(\text{ mod }\pi_m^{p^{m-l}+1}),$$ which, combined with the
equality $E(\pi_l)=E(\pi_m)^{p^{m-l}}$, implies that
$\pi_l\equiv\pi_m^{p^{m-l}} (\text{ mod }\pi_m^{p^{m-l}+1})$.
\begin{corollary}We have
$\pi_{m-i}^{p^j}\equiv\pi_m^{p^{i+j}} (\text{ mod
}\pi_m^{p^{i+j}+1})$.\end{corollary}
\begin{lemma}If $j\leq m-i-1$ and $l<j$, we have
have
$$\text{ord}_p(\frac{\pi_{m-i}^{p^l}}{p^l})
>\text{ord}_p(\frac{\pi_{m-i}^{p^j}}{p^j}).$$\end{lemma}
\begin{corollary}If $j\leq m-i-1$, we have
$p^j\gamma_{i,j}\equiv\pi_m^{p^{i+j}} (\text{ mod
}\pi_m^{p^{i+j}+1})$.\end{corollary}
\begin{corollary}Suppowse that $j\leq m-i-1$. Then
$\text{ord}_p(p^j\gamma_{i,j})>\frac{\deg(p^ju)}{p-1}$ if
$\deg(p^{m-i-1}u)\leq1$, and
$\text{ord}_p(p^j\gamma_{i,j}-\pi^{D\deg(p^ju)})>\frac{\deg(p^ju)}{p-1}$
    if $\deg(p^{m-i-1}u)=1.$\end{corollary}
\begin{lemma}If $j\geq
m-i$, we have
$$\text{ord}_p(p^j\gamma_{i,j})-\frac{\deg(p^ju)}{p-1}\geq
p^{j-(m-i)+1}-1.$$\end{lemma} {\it Proof. }Since
$\gamma_{i,j}=-\sum_{l=j+1}^{\infty}\frac{\pi_{m-i}^{p^{l}}}{p^l}$,
and $\text{ord}_p(\frac{\pi_{m-i}^{p^{l}}}{p^l})
\geq\frac{p^{j+1}}{p^{m-i-1}(p-1)}-{j+1}$ when $j\geq m-i$ and
$l\geq j+1$, we have we have
$\text{ord}_p(p^j\gamma_{i,j})\geq\frac{p^{j+1}}{p^{m-i-1}(p-1)}-1$
if $j\geq m-i$. The lemma now follows from the fact that
$\deg(p^{m-i-1}u)\leq1$.

Write
$$B=\{\sum_{u\in L(\Delta)}a_ux^u\in L(\frac{1}{p-1}):
0\leq \text{ord}_p(a_u)-\frac{\deg(u)}{p-1}\rightarrow
+\infty\text{ as }\deg(u) \rightarrow \infty \}.$$
\begin{corollary}For $k=1,\cdots,n$, we have $\widehat{_kf}\in
B$, and
$$\widehat{_kf}\equiv
\sum\limits_{i=0}^{m-1}\sum\limits_{j=0}^{m-i-1}
\sum\limits_{\deg(p^{m-i-1}u)=1}u_k\omega(a_{iu}^{p^j})\pi^{D\deg(p^ju)}x^{p^ju}\
(\text{mod }\pi B).$$
\end{corollary}
\section{The $p$-adic trace formula}
We shall relate the $L$-function $L_f(t)$ to the characteristic
polynomials of an operator $(p^nF^{-1})^a$ on $p$-adic spaces.

Since $E_f(x)\in L(\frac{1}{p-1})$ (Lemma 3.1), and
$\psi_p:\sum_{u\in L(\Delta)}a_ux^u\mapsto\sum_{u\in
L(\Delta)}a_{pu}x^u$ maps $L(b)$ to $L(pb)$, we have the following
lemma.
\begin{lemma}The map
$p^nF^{-1}: g\mapsto \sigma^{-1}\circ\psi_p(E_f(x)g)$ sends
$L(\frac{1}{p-1})$ to $L(\frac{p}{p-1})$. In particular,
$p^nF^{-1}$ acts on $B$.
\end{lemma}
Note that $p^nF^{-1}$ is $\sigma^{-1}$-linear, and
$(p^nF^{-1})^a=\psi_p^a\circ
\prod\limits_{i=0}^{a-1}E_f^{\sigma^i}(x^{p^i})$ is
$\mathbb{Z}_p[\mu_{q-1},\pi_m,\pi]$-linear. Write
$$
\prod\limits_{i=0}^{ak-1}E_f^{\sigma^i}(x^{p^i})=\sum\limits_{u
\in L(\Delta)}a_ux^u.
$$Then the trace of $(p^nF^{-1})^{ak}$ on $B$ is $\sum\limits_{u\in L(\Delta)}a_{(q^k-1)u}$.
And $$S_k(f)=(-1)^{n-1}(q^k-1)^n\sum\limits_{u \in
L(\Delta)}a_{(q^k-1)u}.$$ So we have the following preliminary
trace formula.
\begin{prop}For
$k=1,2,\cdots$, we have
$$S_k(f)
=-(1-q^k)^n\text{Tr}((p^nF^{-1})^{ak};B).$$Equivalently,
$$L_f(t)=
=\prod\limits_{i=0}^n\det(1-(p^nF^{-1})^aq^it;B)^{(-1)^i(\begin{array}{c}
  n \\
  i \\
\end{array})}$$
\end{prop}Let
$e_1=(1,0,\cdots,0)$, $\cdots$, $e_n=(0,\cdots,0,1)$. For
$l=0,1,\cdots,n$, we write
$$K_l=\bigoplus\limits_{1\leq
i_1<\cdots<i_l\leq n} Be_{i_1}\wedge\cdots\wedge e_{i_l}$$ and
define
$$p^nF^{-1}:K_l\rightarrow K_l,
\ ge_{i_1}\wedge\cdots\wedge e_{i_l}\mapsto
p^{l+n}F^{-1}(g)e_{i_1}\wedge\cdots\wedge e_{i_l}.$$Then the
preliminary trace formula takes the following form.
\begin{prop}For $k=1,2,\cdots$, we have
$$S_k(f)
=\sum\limits_{l=0}^n (-1)^{l+1} \text{Tr}((p^nF^{-1})^{ak};K_l).$$
\end{prop} By
Corollary 3.9, $\hat{D}_j:g\mapsto (x_j\frac{\partial }{\partial
x_j}+\widehat{_jf})g$, $j=1,\cdots,n$, operate on $B$. Obviously,
they commute with each other. So, for $l=1,\cdots,n$,
$$\hat{\partial}:K_l\rightarrow K_{l-1},\
ge_{i_1}\wedge\cdots\wedge
e_{i_l}\mapsto\sum\limits_{k=1}^l(-1)^{k-1}\hat{D}_{i_k}(g)e_{i_1}\wedge\cdots\wedge
\hat{e}_{i_k}\wedge\cdots\wedge e_{i_l},\ i_1<\cdots<i_l$$ are
well-defined, and satisfiy $\hat{\partial}^2=0$. Thus we get a
complex
$$K_n\stackrel{\hat{\partial}}{\rightarrow}K_{n-1}
\stackrel{\hat{\partial}}{\rightarrow}\cdots\stackrel{\hat{\partial}}{\rightarrow}K_0.$$
It is easy to check that
$p^nF^{-1}\circ\hat{\partial}=\hat{\partial}\circ p^nF^{-1}$. That
is, $p^nF^{-1}$ operates on the complex
$(K_{\bullet},\hat{\partial})$. Therefore we have the following
homological trace formula.
\begin{prop}For $k=1,2,\cdots$, we have
$$S_k(f)=\sum\limits_{l=0}^n (-1)^{l+1}
\text{Tr}((p^nF^{-1})^{ak};H_l(K_{\bullet},\hat{\partial})).$$
Equivalently, $$L_f(t)=\prod\limits_{l=0}^n
\det(1-(p^nF^{-1})^at;H_l(K_{\bullet},\hat{\partial}))^{(-1)^l}.$$
\end{prop}
\section{The total degree of the $L$-function}
We shall study the Newton polygon of $\det(1-(p^nF^{-1})^at;B)$,
and then prove Theorem 1.2.
\begin{prop}The Newton polygon of
$\det(1-(p^nF^{-1})^at;B)$ with respect to $\text{ord}_q$ lies
above the Hodge polygon of
$\sum\limits_{i=0}^{+\infty}W_{\Delta}(k)t^k$ of degree
$D$.\end{prop} Write $E_f(x)=\sum\limits_{u\in
L(\Delta)}a_u\pi^{D\deg(u)}x^u$,
$a_u\in\mathbb{Z}_p[\mu_{q-1},\pi_m,\pi]$. Then the matrix of
$p^nF^{-1}$ with respect to the orthonormal basis
$\{\pi^{D\deg(u)}x^u\}_{u\in L(\Delta)}$, written as a column
vector, is
$$A^{\sigma^{-1}}=(a_{pw-u}^{\sigma^{-1}}\pi^{D((p-1)\deg(w)+c(w,u))})_{w,u},\
c(w,u)=\deg(pw-u)+\deg(u)-p\deg(w)\geq 0.$$ So, the matrix of
$(p^nF^{-1})^a$ with respect to that orthonormal basis is
$AA^{\sigma}\cdots A^{\sigma^{a-1}}$. Obviously, the Newton
polygon of $\det(1-At)$ with respect to $\text{ord}_p$ lies above
the polygon with vertices at points $(0,0)$ and
$$(\sum\limits_{i=0}^kW_{\Delta}(i),
\sum\limits_{i=0}^kW_{\Delta}(i)\frac{i}{D}),\ k=0,1,\cdots.$$ It
follows that the Newton polygon of
$\det(1-(p^nF^{-1})^at;B)=\det(1-AA^{\sigma}\cdots
A^{\sigma^{a-1}}t)$ with respect to $\text{ord}_q$ lies above the
polygon with vertices at points $(0,0)$ and
$$(\sum\limits_{i=0}^kW_{\Delta}(i),
\sum\limits_{i=0}^kW_{\Delta}(i)\frac{i}{D}),\ k=0,1,\cdots.$$ The
proposition is proved.
\begin{corollary}If $j\leq n+1$, then
$\det(1-(p^nF^{-1})^at;B)$ has at most
$\sum\limits_{k=0}^{Dj}W_{\Delta}(k)$ zeros of $q$-order $\leq
j-1$.\end{corollary}{\it Proof. }Define
$$\sum\limits_{k=0}^{+\infty}h_{\Delta}(k)t^k=
(1-t)^n\sum\limits_{k=0}^{+\infty}W_{\Delta}(k)t^k.$$ Since
$\sum\limits_{k=0}^{+\infty}h_{\Delta}(k)t^k$ is a polynomial of
degree $\leq n$ with nonnegative coefficients by a lemma of
Kouchnirenko [Ko, Lemma 2.9], and
$$\sum\limits_{k=0}^{jD-i}(\begin{array}{c}
  n-1+k \\
  n-1 \\
\end{array})(k+i)=(\begin{array}{c}
  n+Dj-i \\
  n \\
\end{array})(\frac{n(Dj-i)}{n+1}+i)\geq(\begin{array}{c}
  n+Dj-i \\
  n \\
\end{array})D(j-1),$$ we have
$$\frac{1}{D}\sum\limits_{k=0}^{jD}kW_{\Delta}(k)
=\frac{1}{D}\sum\limits_{k=0}^{jD}k\sum\limits_{i=0}^kh_{\Delta}(i)(\begin{array}{c}
  n-1+k-i \\
  n-1 \\
\end{array})$$
$$=\frac{1}{D}\sum\limits_{i=0}^nh_{\Delta}(i)\sum\limits_{k=i}^{jD}(\begin{array}{c}
  n-1+k-i \\
  n-1 \\
\end{array})k=\frac{1}{D}\sum\limits_{i=0}^nh_{\Delta}(i)\sum\limits_{k=0}^{jD-i}(\begin{array}{c}
  n-1+k \\
  n-1 \\
\end{array})(k+i)$$
$$\geq(j-1)\sum\limits_{i=0}^nh_{\Delta}(i)\sum\limits_{k=0}^{jD-i}(\begin{array}{c}
  n-1+k \\
  n-1 \\
\end{array})\geq(j-1)\sum\limits_{k=0}^{jD}W_{\Delta}(k).$$
The corollary now follows from the above inequality by Proposition
5.1.

We now prove Theorem 1.2. Since the reciprocal zeros and
reciprocal poles of $L_f(t)$ are of $q$-order $\leq n$, its total
number, by the preliminary trace formula, is bounded by the number
of reciprocal zeros of
$\prod\limits_{i=0}^n\det(1-(p^nF^{-1})^aq^it;B)^{(\begin{array}{c}
  n \\
  i \\
\end{array})}$. By Corollary 5.2, that number is bounded by
$$\sum\limits_{i=0}^n(\begin{array}{c}
  n \\
  i \\
\end{array})\sum\limits_{k=0}^{D(n-i+1)}W_{\Delta}(k).$$ Theorem
1.2 is proved.
\section{The acyclicity of the $p$-adic complex}
In this section we shall prove the following proposition, which
implies the first statement of Theorem 1.3.
\begin{prop}If $f$ non-degenerate with respect to
$\Delta_{\infty}(f)$, then $(K_{\bullet},\hat{\partial})$ is
acyclic at positive dimensions, and
$H_0(K_{\bullet},\hat{\partial})$ is a
$\mathbb{Z}_p[\mu_{q-1},\pi_m,\pi]$-module free of rank
$n!\text{Vol}(\Delta_{\infty}(f))$.\end{prop}

Write
$$\bar{B}:=\mathbb{F}_q[x^{L(\Delta)}]:=\{\sum\limits_{u\in L(\Delta)}a_ux^u:a_u\in\mathbb{F}_q\}.$$
It is a ring with the multiplication rule
$$x^ux^{u'}=\left\{\begin{array}{ll}
x^{u+u'},&\text{ if } u \text{ and }u'\text{ are cofacial}, \\
0, &\text{ otherwise.}
\end{array}\right. $$
Define $$B\rightarrow\bar{B},\ \sum\limits_{u\in
L(\Delta)}a_u\pi^{D\deg(u)}x^u\mapsto\sum\limits_{u \in
L(\Delta)}\bar{a}_ux^u,$$ where $\bar{a}_u$ is the residue class
of $a_u$ modulo the maximal ideal of
$\mathbb{Z}_p[\mu_{q-1},\pi_m,\pi]$.
\begin{lemma}The map $B \rightarrow \bar{B}$ is a ring homomorphism.
And the
sequence
$$0\rightarrow B\rightarrow B
\rightarrow\bar{B}\rightarrow0$$ is exact.\end{lemma} For
$j=1,\cdots,n$, we define
$$\bar{D}_j: \bar{B}\rightarrow\bar{B},\
g\mapsto(x_j\frac{\partial }{\partial x_j}+\overline{_jf})g,$$
where
$$\overline{_jf}=\sum\limits_{i=0}^{m-1}\sum\limits_{j=0}^{m-i-1}
\sum\limits_{\deg(p^{m-i-1}u)=1}u_k a_{iu}^{p^j}x^{p^ju}.$$ By
Corollary 3.8, we have the following lemma. \begin{lemma}For
$j=1,\cdots,n$, the diagram
$$\begin{array}{ccc}
  B &\rightarrow& \bar{B} \\
 \hat{D}_j\downarrow &  & \bar{D}_j\downarrow \\
  B & \rightarrow & \bar{B} \\
\end{array}$$
is commutative.\end{lemma}
For $l=0,\cdots,n$, we define
$$\bar{K}_l=\bigoplus\limits_{1\leq i_1<\cdots<i_l\leq n}
\bar{B}e_{i_1}\wedge\cdots\wedge e_{i_l}.$$
For $l=1,\cdots,n$,
we define
$$\bar{\partial}:\bar{K}_l\rightarrow\bar{K}_{l-1},\ ge_{i_1}\wedge\cdots\wedge e_{i_l}
\mapsto \sum\limits_{k=1}^l(-1)^{k-1}\bar{D}_{i_k}(g)
e_{i_1}\wedge\cdots\wedge\hat{e}_{i_k}\wedge\cdots\wedge e_{i_l},\
i_1<\cdots<i_l.$$ It is easy to see that the sequence
$$\bar{K}_n\stackrel{\bar{\partial}}{\rightarrow}\bar{K}_{n-1}
\stackrel{\bar{\partial}}{\rightarrow}\cdots\stackrel{\bar{\partial}}{\rightarrow}\bar{K}_0$$
is a complex.
\begin{prop}The map $B\rightarrow\bar{B}$ induces a morphism of
complexes from $(K_{\bullet},\hat{\partial})$ to
$(\bar{K}_{\bullet},\bar{\partial})$. Moreover, the sequence
$$0\rightarrow (K_{\bullet},\hat{\partial}) \rightarrow
(K_{\bullet},\hat{\partial})  \rightarrow
(\bar{K}_{\bullet},\bar{\partial})\rightarrow0$$ is
exact.\end{prop} {\it Proof. }The first statement follows from
Lemma 6.3, and the second follows from Lemma 6.2.

By Proposition 6.4, and a lemma of Monsky [Mo, Theorem 8.5], the
proof of Proposition 6.1 is reduced to the proof of the following
proposition.
\begin{prop}If $f$ is non-degenerate with respect to
$\Delta_{\infty}(f)$, then $(\bar{K}_{\bullet},\bar{\partial})$ is
acyclic at positive dimensions, and
$H_0((\bar{K}_{\bullet},\bar{\partial}))$ is a
$\mathbb{F}_q$-vector space of dimension
$n!\text{Vol}(\Delta_{\infty}(f))$.\end{prop}
For $j=1,\cdots,n$,
we define
$$\overline{_jf}^0=\sum\limits_{i=0}^{m-1}
\sum\limits_{\deg(p^{m-i-1}u)=1}u_k
a_{iu}^{p^{m-i-1}}x^{p^{m-i-1}u}.$$ For $l=1,\cdots,n$, we
define$$\bar{\partial}^0:\bar{K}_l\rightarrow\bar{K}_{l-1},\
ge_{i_1}\wedge\cdots\wedge e_{i_l}\mapsto
\sum\limits_{k=1}^l(-1)^{k-1}\overline{_{i_k}f}^0g
e_{i_1}\wedge\cdots\wedge\hat{e}_{i_k}\wedge\cdots\wedge e_{i_l},\
i_1<\cdots<i_l.$$ Then
$$\bar{K}_n\stackrel{\bar{\partial}^0}{\rightarrow}\bar{K}_{n-1}
\stackrel{\bar{\partial}^0}{\rightarrow}
\cdots\stackrel{\bar{\partial}^0}{\rightarrow}\bar{K}_0$$ is a
complex. In the next section, we shall prove the following
proposition.
\begin{prop}If $f$ is non-degenerate with respect to
$\Delta:=\Delta_{\infty}(f)$, then the complex
$(\bar{K}_{\bullet},\bar{\partial}^0)$ is acyclic at positive
dimensions, and the Poincar\'{e} series of
$H_0((\bar{K}_{\bullet},\bar{\partial}^0))$ is $P_{\Delta}(t)$. In
particular, $H_0((\bar{K}_{\bullet},\bar{\partial}^0))$ is a
$\mathbb{F}_q$-vector space of dimension $n!\text{Vol}(\Delta)$.
\end{prop}

We now deduce the first statement of Proposition 6.5 from
Proposition 6.6. In a given
 a homology class of positive dimension,
we choose one representative $\xi$ of lowest degree. We claim that
$\xi=0$. Otherwise, let $\xi^0$ be the leading term of $\xi$. We
have $\bar{\partial}^0(\xi^0)=0$ since it is the leading term of
$\bar{\partial}(\xi)=0$. By the acyclicity of
$(\bar{K}_{\bullet},\bar{\partial}^0)$,
$\xi^0=\bar{\partial}^0(\eta)$ for some $\eta$. The form
$\xi-\bar{\partial}(\eta)$ is now of lower degree than $\xi$,
contradicting to our choice of $\xi$. The proposition is proved.

The second statement of Proposition 6.5 follows the following
proposition.
\begin{prop}Let $V$ be a basis of $\bar{K}_0$ modulo
$\bar{\partial}^0(\bar{K}_1)$ consisting of homogeneous elements.
Then $V$ is also a basis of $\bar{K}_0$ modulo
$\bar{\partial}(\bar{K}_1)$.
\end{prop}
{\it Proof. }First, we show that $\bar{K}_0$ is generated by $V$
and $\bar{\partial}(\bar{K}_1)$. Otherwise, among elements of
$\bar{K}_0$ which are not linear combinations of elements of $V$
and $\bar{\partial}(\bar{K}_1)$, we choose one of lowest degree.
We may suppose that it is of form $\bar{\partial}^0(\xi)$. Let
$\xi^0$ be the leading term of $\xi$. Then
$\bar{\partial}^0(\xi)-\bar{\partial}(\xi^0)$ is not a linear
combination of elements of $V$ and $\bar{\partial}(\bar{K}_1)$,
and is of lower degree than $\partial^0(\xi)$. This is a
contradiction. Therefore $\bar{K}_0$ is generated by $E$ and
$\bar{\partial}(\bar{K}_1)$. It remains to show that $\xi=0$
whenever $\xi$ belongs to $\bar{\partial}(\bar{K}_1)$ and is a
linear combination of elements of $V$. Otherwise, we may choose
one element $\zeta$ of lowest degree such that
$\xi=\bar{\partial}(\zeta)$. Let $\zeta^0$ be the leading term of
$\zeta$. Then $\bar{\partial}^0(\zeta^0)$ is a linear combination
of elements of $V$ since it is the leading term of
$\bar{\partial}(\zeta)$. So we have $\bar{\partial}^0(\zeta^0)=0$.
By the acyclicity of $(\bar{K}_{\bullet},\bar{\partial}^0)$,
$\zeta^0=\bar{\partial}^0(\eta)$ for some $\eta$. The form
$\zeta-\bar{\partial}(\eta)$ is now of lower degree than $\zeta$,
contradicting to our choice of $\zeta$. This completes the proof
of the proposition.
\section{The complex obtained by reduction}
In this section, we shall prove Proposition 6.6. The second
statement follows from the first, and the last follows from the
second and a lemma of Kouchnirenko [Ko, Lemma 2.9]. So it remains
to prove the acyclicity of the complex
$(\bar{K}_{\bullet},\bar{\partial}^0)$.

Let $\tau$ be a face of $\Delta$ that does not contain the origin,
and $\bar{\tau}$ is the convex hull in $\mathbb{Q}^n$ generated by
$\tau$ and the origin. For $\alpha_1,\cdots,\alpha_s$ in
$\mathbb{F}_q[x^{L(\bar{\tau})}]$, we define
$\bar{K}_{\bullet}(\bar{\tau},\{\alpha_j\}_{j=1}^s)$ to be the
complex$$\bar{K}_l(\bar{\tau},\{\alpha_j\}_{j=1}^s)=\bigoplus\limits_{1\leq
i_1<\cdots<i_l\leq
s}\mathbb{F}_q[x^{L(\bar{\tau})}]e_{i_1}\wedge\cdots\wedge
e_{i_l},\ l=0,\cdots,s$$ with derivation
$$ge_{i_1}\wedge\cdots\wedge e_{i_l}\mapsto
\sum\limits_{k=1}^l(-1)^{k-1}\alpha_{i_k}g
e_{i_1}\wedge\cdots\wedge\hat{e}_{i_k}\wedge\cdots\wedge e_{i_l},\
1\leq i_1<\cdots<i_l\leq s.$$ By a proposition of Kouchnirenko
[Ko, Proposition 2.6] and the argument of Adolphson-Sperber [AS2,
p379], the sequence
$$0\rightarrow\bar{K}_{\bullet}^{\emptyset}(f)\rightarrow
\bigoplus\limits_{\dim\tau=n-1}\bar{K}_{\bullet}(\bar{\tau},
\{\overline{_jf}^{\tau}\}_{j=1}^n)\rightarrow\cdots\rightarrow
\bigoplus\limits_{\dim\tau=0}\bar{K}_{\bullet}(\bar{\tau},
\{\overline{_jf}^{\tau}\}_{j=1}^n)\rightarrow
\bar{K}_{\bullet}^{-1}\rightarrow0$$ is exact, where $\tau$
denotes a face of $\Delta$ that does not contain the origin, and
$$\bar{K}_l^{-1}=\left\{%
\begin{array}{ll}
    \mathbb{F}_q^{(\begin{array}{c}
  n \\
  l \\
\end{array})}, & \text{ if the origin is in the interior of }\Delta\text{ and }1\leq l\leq n, \\
    0, & \text{otherwise.} \\
\end{array}%
\right. $$
By the exactness of that sequence, the acyclicity of the
complex $(\bar{K}_{\bullet},\bar{\partial}^0)$ follows from the
following lemma.
\begin{lemma}Let $f$ be a Witt vector of length
$m$ with coefficients in
$\mathbb{F}_q[x_1^{\pm1},\cdots,x_n^{\pm1}]$. Suppose that $f$ is
non-degenerate with respect to $\Delta:=\Delta_{\infty}(f)$ and
$\dim\Delta=n$. Let $\tau$ be a face of $\Delta$ of dimension
$s-1$ that does not contain the origin. Then the complex
$\bar{K}_{\bullet}(\bar{\tau}, \{\overline{_jf}^{\tau}\}_{j=1}^n)$
is acyclic at dimensions $>n-s$.
\end{lemma}
Since the
sequence
$$0\rightarrow\bar{K}_{\bullet}(\bar{\tau},
\{\alpha_j\}_{j=1}^{s-1})\rightarrow\bar{K}_{\bullet}(\bar{\tau},
\{\alpha_j\}_{j=1}^s)\rightarrow \bar{K}_{\bullet}(\bar{\tau},
\{\alpha_j\}_{j=1}^{s-1})[-1]\rightarrow0$$ is exact, Lemma 7.1
follows from the following one.
\begin{lemma}Suppose that $f$ is
non-degenerate with respect to $\Delta_{\infty}(f)$. If $\tau$ is
a face of $\Delta$ of dimension $r-1$ that does not contain the
origin, then there are $1\leq i_1<\cdots<i_r\leq n$ such that the
complex
$\bar{K}_{\bullet}(\bar{\tau},\{\overline{_{i_j}f}^{\tau}\}_{j=1}^r)$
is acyclic at positive dimensions.
\end{lemma}
{\it Proof. }There are $\{i_1,\ldots,i_r\}\subset\{1,\cdots,n\}$
and $(\alpha_{kj})\in\mathbb{Q}\cap\mathbb{Z}_p$ ($1\leq k\leq n$,
$1\leq j\leq r$) such that
$u_k=\alpha_{1k}u_{i_1}+\cdots+\alpha_{rk}u_{i_r}$ for all
$u=(u_1,\ldots,u_n)\in L(\bar{\tau})$. Let $\sigma$ be any face of
$\tau$. We have
$$\overline{_jf}^{\sigma}=\alpha_{1j}
\overline{_{i_1}f}^{\sigma}+\cdots+\alpha_{rj}\overline{_{i_r}f}^{\sigma}.$$
So $\overline{_{i_1}f}^{\sigma}$, $\cdots$,
$\overline{_{i_r}f}^{\sigma}$ have no common zeros in
$(\overline{\mathbb{F}}_q^{\times})^n$. By a theorem of
Kouchnirenko [Ko, Theorem 6.2], $\overline{_{i_1}f}^{\tau}$,
$\cdots$, $\overline{_{i_r}f}^{\tau}$  generate in
$\mathbb{F}_q[x^{L(\bar{\tau})}]$ an ideal of finite codimension.
Note that $\mathbb{F}_q[x^{L(\bar{\tau})}]$ is Cohen-Macaulay by a
theorem of Hochster [Ho, Theorem 1]. The complex
$\bar{K}_{\bullet}(\bar{\tau},\{\overline{_{i_j}f}^{\tau}\}_{j=1}^r)$
is acyclic at positive dimensions by a theorem of Serre [Se,
Theorem 3, Chapter IV]. The lemma is proved.
\section{The Newton polygon of the $L$-function}In this section we
shall prove the second statement of Theorem 1.3. (The last
statement follows from the second by a lemma of Kouchnirenko [Ko,
Lemma 2.9].) By the argument of Dwork [Dw2, \S7], it suffices to
prove the following proposition.
\begin{prop}If $f$
is non-degenerate with respect to $\Delta:=\Delta_{\infty}(f)$,
then the Newton polygon of
$\det(1-p^nF^{-1}t;H_0(K_{\bullet},\hat{\partial}))$ with respect
to $\text{ord}_p$ lies above the Hodge polygon of $P_{\Delta}(t)$
of degree $D(\Delta)$, and their endpoints coincide.
\end{prop}

Let $\bar{V}$ be a basis of $\bar{K}_0$ modulo
$\bar{\partial}^0(K_1)$ consisting of homogeneous elements. By
Proposition 6.7, it is also a basis of $\bar{K}_0$ modulo
$\bar{\partial}(K_1)$. Define
$$V=\{\sum \omega(a_u)x^u:\sum a_ux^u\in\bar{V}\}.$$It is a basis of $B$ modulo
$\sum\limits_{k=1}^n\hat{D}_kB$. For real numbers
$b>\frac{1}{p-1}$ and $c$, we write
$$L(b,c)=\{\sum_{u\in L(\Delta)}a_ux^u :
 a_u\in \mathbb{Q}_p[\mu_{q-1},\pi_m,\pi], \text{ ord}_p(a_u)\geq b\deg(u)+c\}.$$
It is compact with respect to the topology of coefficientwise
convergence. Let $V(b,c)$ be the subset of elements of $L(b,c)$
which are finite linear combinations of elements of $V$. In the
next section we shall prove the following proposition.
\begin{prop}If $\frac{1}{p-1}<b<\frac{p}{p-1}$, then
$$L(b,c)=V(b,c)+\sum\limits_{k=1}^n\hat{D}_kL(b,c+b-\frac{1}{p-1}).$$
\end{prop}

We now prove the first statement of Proposition 8.1. For each
$\xi\in V$, we write
$$p^nF^{-1}(\pi^{D\deg(\xi)}\xi)\equiv\sum\limits_{\eta\in
V}c_{\eta,\xi}\pi^{D\deg(\eta)}\eta \ (\text{mod
}\sum\limits_{k=1}^n\hat{D}_kB),\
c_{\eta,\xi}\in\mathbb{Z}_p[\mu_{q-1},\pi_m,\pi]$$ By Lemma 4.1,
$p^nF^{-1}(\pi^{D\deg(\xi)}\xi)$ lies in the space
$L(\frac{p}{p-1})$. So, by Proposition 8.2,
$c_{\eta,\xi}\pi^{D\deg(\eta)}\eta$ lies in every $L(b)$ with
$\frac{1}{p-1}<b<\frac{p}{p-1}$. That is,
$\text{ord}_p(c_{\eta,\xi})\geq(b-\frac{1}{p-1})\deg(\eta)$ for
every $\frac{1}{p-1}<b<\frac{p}{p-1}$. Thus we have
$\text{ord}_p(c_{\eta,\xi})\geq\deg(\eta)$. Therefore, the Newton
polygon of the characteristic polynomial of $(c_{\eta,\xi})$,
which is now the Newton polygon of
$\det(1-p^nF^{-1}t;H_0(K_{\bullet},\hat{\partial}))$, lies above
the Hodge polygon of $P_{\Delta}(t)$ of degree $D$. In particular,
$\text{ord}_p(\det(c_{\eta,\xi}))\geq\sum\limits_{\xi\in
V}\deg(\xi)$.

It remains to show that the Newton polygon of
$\det(1-p^nF^{-1}t;H_0(K_{\bullet},\hat{\partial}))$ share the
same endpoints with the Hodge polygon of $P_{\Delta}(t)$ of degree
$D$. Define
$$\phi_p:L(\frac{p}{p-1})\rightarrow L(\frac{1}{p-1}),\
\sum\limits_{u\in L(\Delta)}a_ux^u\mapsto\sum\limits_{u\in
L(\Delta)}a_ux^{pu}.$$ Obviuosly,
$p^nF^{-1}\circ(E_f^{-1}\circ\phi_p\circ\sigma)=1$ on
$L(\frac{p}{p-1})$. At the end of this we shall prove the
following proposition.
\begin{prop}Modulo
$\sum\limits_{k=1}^n\hat{D}_kL(\frac{1}{p-1})$,
 the space $L(\frac{1}{p-1})$ is generated by
$\{\pi^{D\deg(\xi)}\xi:\xi\in V\}$.
\end{prop}
So, for each $\xi\in V$, we can find
$b_{\eta,\xi}\in\mathbb{Z}_p[\mu_{q-1},\pi_m,\pi]$
$$E_f^{-1}\circ\phi_p\circ\sigma(\pi^{D\deg(\xi)p}\xi)\equiv\sum\limits_{\eta\in
V}b_{\eta,\xi}\pi^{D\deg(\eta)}\eta \ (\text{mod
}\sum\limits_{k=1}^n\hat{D}_kL(\frac{1}{p-1})).$$ It follows that
$(c_{\eta,\xi})(b_{\eta,\xi})=\text{diag}\{\pi^{D\deg(\xi)(p-1)},\xi\in
V\}$. So $\text{ord}_p(\det(c_{\eta,\xi}))\leq\sum\limits_{\eta\in
V}\deg(\eta)$. Therefore
$$\text{ord}_p(\det(c_{\eta,\xi}))=\sum\limits_{\eta\in
V}\deg(\eta).$$ That is, the Newton polygon of
$\det(1-p^nF^{-1}t;H_0(K_{\bullet},\hat{\partial}))$ share the
same endpoints with the Hodge polygon of $P_{\Delta}(t)$ of degree
$D$.

We now prove Proposition 8.3. Let
$\xi=\sum\limits_{u\in L(\Delta)}a_ux^u\in L(\frac{p}{p-1})$. For
$N=0,1,\cdots$, write $\xi^{(N)}=\sum\limits_{u\in
L(\Delta),\deg(u)\leq N}a_ux^u\in B$. As
$\{\pi^{D\deg(\eta)}\eta:\eta\in V\}$ is a basis of $B$ modulo
$\sum\limits_{k=1}^n\hat{D}_kB$, there are elements
$\xi_k^{(N)}\in B$ ($k=1,\cdots,n$) such that
$$\xi^{(N)}-\sum\limits_{k=1}^n\hat{D}_k\xi_k^{(N)}=\sum\limits_{\eta\in
V}a_{\eta}^{(N)}\pi^{D\deg(\eta)}\eta.$$ As $L(\frac{1}{p-1})$ is
compact with respect to the topology of coefficientwise
convergence, the sequence
$(\{\xi_k^{(N)}\}_{k=1}^n,\{a_{\eta}^{(N)}\}_{\eta\in V})$ ,
$N=0,1,\cdots$, has an adherent point
$(\{\xi_k\}_{k=1}^n,\{a_{\eta}\}_{\eta\in V})$ in the space
$L(b)^n\times(\mathbb{Z}_p[\mu_{q-1},\pi_m,\pi])^{|V|}$. Therefore
we get
$$\xi-\sum\limits_{k=1}^n\hat{D}_k\xi_k=\sum\limits_{\eta\in
V}a_{\eta}\pi^{D\deg(\eta)p^{m-1}}\eta.$$ This completes the proof
of Proposition 8.3.
\section{The space $L(b,c)$}
In this section we shall prove Propositions 8.2.

For $k=1,\cdots,n$, we write
$$\widehat{_kf}^0=\sum\limits_{i=0}^{m-1}p^{m-i-1}\gamma_{i,m-i-1}
\sum\limits_{\deg(p^{m-i-1}u)=1}u_k\omega(a_{iu}^{p^{m-i-1}})
x^{p^{m-i-1}u}.$$
For $l=1,\cdots,n$, we define
$$\hat{\partial}^0:K_l\rightarrow K_{l-1},\
ge_{i_1}\wedge\cdots\wedge
e_{i_l}\mapsto\sum\limits_{k=1}^l(-1)^{k-1}\widehat{_{i_k}f}^0ge_{i_1}\wedge\cdots\wedge
\hat{e}_{i_k}\wedge\cdots\wedge e_{i_l},\ i_1<\cdots<i_l.$$ It is
easy to see that
$$0\rightarrow
(K_{\bullet},\hat{\partial}^0) \rightarrow
(K_{\bullet},\hat{\partial}^0) \rightarrow
(\bar{K}_{\bullet},\bar{\partial})\rightarrow0$$ is an exact
sequence of complexex. So we have the following lemma.
\begin{lemma}Modulo
$\sum\limits_{k=1}^n\widehat{_kf}^0B$, the space $B$ is generated
by $\{\pi^{D\deg(\xi)}\xi:\xi\in V\}$.
\end{lemma}
\begin{corollary}If $b>\frac{1}{p-1}$, then
$$L(b,c)=V(b,c)+\sum\limits_{k=1}^n\widehat{_kf}^0L(b,c+b-\frac{1}{p-1}).$$
\end{corollary}
{\it Proof. }Let $\xi\in L(b,c)$, $\xi_v$ ($v\in \deg L(\Delta)$)
its homogeneous part of degree $v$, and $k_v$ the least integer
such that $\text{ord}_p(\pi^{k_v})\geq bv+c$. Then
$\pi^{Dv-k_v}\xi_v\in B$. By the above lemma, we may write
$$\pi^{Dv-k_v}\xi_v=\sum\limits_{\eta\in V,\deg(\eta)\leq
v}a_{\eta}^{(v)}
\pi^{D\deg(\eta)}\eta+\sum\limits_{i=1}^n\widehat{_if}^0\eta_i^{(v)},$$
where $a_{\eta}^{(v)}\in \mathbb{Z}_p[\mu_{q-1},\pi_m,\pi]$, and
$\eta_i^{(v)}\in B$ is of degree $\leq v-1$. It follows that
$$\xi=\sum\limits_{\eta\in V}\eta\pi^{D\deg(\eta)}\sum\limits_{v\geq\deg(\eta)}a_{\eta}^{(v)}
\pi^{k_v-Dv}+
\sum\limits_{i=1}^n\widehat{_if}^0\sum\limits_{v\in\deg
L(\Delta)}\pi^{k_v-Dv}\eta_i^{(v)}.$$ It is easy to see that the
first term on the right-hand side converges to an element in
$V(b,c)$, and the inner sum in the second term converges to an
element in $L(b,c+b-\frac{1}{p-1})$. The corollary is proved.

For $k=1,\cdots,n$, we define $$D_k: B\rightarrow B,\
g\mapsto(x_k\frac{\partial }{\partial
x_k}+\sum\limits_{i=0}^{m-1}\sum\limits_{j=0}^{m-i-1}
p^j\gamma_{i,j}\sum\limits_{u\in I_i}u_k
\omega(a_{iu}^{p^j})x^{p^ju})g.$$ For $l=1,\cdots,n$, we define
$$\partial:K_l\rightarrow K_{l-1},\
ge_{i_1}\wedge\cdots\wedge
e_{i_l}\mapsto\sum\limits_{k=1}^l(-1)^{k-1}D_k(g)e_{i_1}\wedge\cdots\wedge
\hat{e}_{i_k}\wedge\cdots\wedge e_{i_l},\ i_1<\cdots<i_l.$$
\begin{corollary}If $b>\frac{1}{p-1}$, then
$$L(b,c)=V(b,c)+\sum\limits_{k=1}^nD_kL(b,c+b-\frac{1}{p-1}).$$
\end{corollary}
{\it Proof. }Note that $\widehat{_kf}^0-D_k$ maps $L(b,c)$ to
$L(b,c-(b-\frac{1}{p-1})e)$ for some constant $e<1$. Let $\xi\in
L(b,c)$. By the previous corollary and induction, we can find a
sequence
$$(\eta_0^{(i)},\cdots,\eta_n^{(i)})\in
V(b,c+i(1-e)(b-\frac{1}{p-1}))\times
L(b,c+(i(1-e)+1)(b-\frac{1}{p-1}))^n,\ i=0,1,\cdots$$ such that
$$\xi=\eta_0^{(0)}+\sum\limits_{k=1}^n\widehat{_kf}^0\eta_k^{(0)},$$
and $$\sum\limits_{k=1}^n(\widehat{_kf}^0-D_k)\eta_k^{(i)}
=\eta_0^{(i+1)}+\sum\limits_{k=1}^n\widehat{_kf}^0\eta_k^{(i+1)}.$$
One sees immediately that $\sum\limits_{i=0}^{\infty}\eta_0^{(i)}$
converges to an element $\eta_0$ in $V(b,c)$, and
$\sum\limits_{i=0}^{\infty}\eta_k^{(i)}$ converges to an element
$\eta_k$ in $L(b,c+b-\frac{1}{p-1})$. Moreover, we have
$\xi=\eta_0+\sum\limits_{k=1}^nD_k\eta_k$. The corollary is
proved.

We now prove Proposition 8.2.
Note that $D_k-\hat{D}_k\in
L(b,p(\frac{p}{p-1}-b)-1)$ by Lemma 3.7. So it maps $L(b,c)$ to
the space $L(b,c+p(\frac{p}{p-1}-b)-1)$. Let $\xi\in L(b,c)$. By
the previous corollary and induction, we can find a sequence
$$(\eta_0^{(i)},\cdots,\eta_n^{(i)})\in
V(b,c+i(p-(p-1)b))\times L(b,c+i(p-(p-1)b)+(b-\frac{1}{p-1}))^n,\
i=0,1,\cdots$$ such that
$$\xi=\eta_0^{(0)}+\sum\limits_{k=1}^nD_k\eta_k^{(0)},$$
and $$\sum\limits_{k=1}^n(D_k-\hat{D}_k)\eta_k^{(i)}
=\eta_0^{(i+1)}+\sum\limits_{k=1}^nD_k\eta_k^{(i+1)}.$$ One sees
immediately that $\sum\limits_{i=0}^{\infty}\eta_0^{(i)}$
converges to an element $\eta_0$ in $V(b,c)$, and
$\sum\limits_{i=0}^{\infty}\eta_k^{(i)}$ converges to an element
$\eta_k$ in $L(b,c+b-\frac{1}{p-1})$. Moreover, we have
$\xi=\eta_0+\sum\limits_{k=1}^nD_k\eta_k$. This completes the
proof of Proposition 8.2.
\section{The weights of the $L$-function}We shall
prove Theorem 1.4.

Let $\Delta$ be a convex polyhedron in $\mathbb{Q}^n$ of dimension
$n$ that contains the origin, and $S_{\Delta}$ the sum of the
volumes of all its $(n-1)$-dimensional faces that contain 0. Write
$$(1-t^{D(\Delta)})^n\sum\limits_{i=0}^{+\infty}W_{\Delta}(i)t^i
=\sum\limits_{i=0}^{D(\Delta)n}h_{\Delta}(i)t^i.$$\begin{lemma}We
have
$$\frac
{1}{D(\Delta)}\sum\limits_{i=0}^{nD(\Delta)}ih_{\Delta}(i)=\frac
n2n!\text{Vol}(\Delta)-\frac{(n-1)!}{2}S_{\Delta}.$$In particular,
$\frac
{1}{D(\Delta)}\sum\limits_{i=0}^{nD(\Delta)}ih_{\Delta}(i)=\frac
n2n!\text{Vol}(\Delta)$ if the origin is an interior point of
$\Delta$.\end{lemma} {\it Proof. }Note that
$$W_{\Delta}(i)=\sum\limits_{k=0}^{i/D(\Delta)}h_{\Delta}(i-D(\Delta)k)\left(%
\begin{array}{c}
  n-1+k \\
  n-1 \\
\end{array}%
\right).$$ So
$$\sum\limits_{i\leq D(\Delta)x}(x-\frac {i}{D(\Delta)})W_{\Delta}(i)
=\sum\limits_{j=0}^{nD(\Delta)}h_{\Delta}(j)
\sum\limits_{0\leq k\leq x-\frac {j}{D(\Delta)}}(x-\frac {j}{D(\Delta)}-k)\left(%
\begin{array}{c}
  n-1+k \\
  n-1 \\
\end{array}%
\right)$$
$$=\sum\limits_{j=0}^{nD(\Delta)}h_{\Delta}(j)
(\frac {x^{n+1}}{(n+1)!}+(\frac n2-\frac {j}{D(\Delta)})\frac
{x^n}{n!}+O(x^{n-1})) .$$ On the other hand, by [AS, (4.12-13)],
$$\sum\limits_{i\leq D(\Delta)x}(x-\frac {i}{D(\Delta)})W_{\Delta}(i)
=n!\text{Vol}(\Delta)\frac{x^{n+1}}{(n+1)!}+\frac{(n-1)!}{2}S_{\Delta}\frac{x^n}{n!}+O(x^{n-1}).$$
The lemma now follows.

We now prove Theorem 1.4. Let $\alpha_i$, $i=1,\cdots,
n!\text{Vol}(\Delta)$, be the eigenvalues of $q^nF^{-1}$ on
$H_0(K_{\bullet},\hat{\partial})$. By Theorem 1.2 and Lemma 10.1,
$$\text{ord}_q(\prod\limits_{i=1}^{n!\text{Vol}(\Delta)}\alpha_i)
=\frac{n}{2}n!\text{Vol}(\Delta).$$ It is known that the
eigenvalues $\alpha_i$ are $l$-adic units when $l$ is a prime
different from $p$. So, by the product formula, we have
$$\prod\limits_{i=1}^{n!\text{Vol}(\Delta)}|\alpha_i|
=q^{\frac{n}{2}n!\text{Vol}(\Delta)}.$$  By a theorem of Kedlaya
[Ke, Theorem 5.6.2], the Frobenius $F$ on
$H_0(K_{\bullet},\hat{\partial})\otimes\mathbb{Q}_p[\mu_{q-1},\pi_m,\pi]$
is of mixed weight $\geq n$. So $q^nF^{-1}$ on
$H_0(K_{\bullet},\hat{\partial})$ is of mixed weight $\leq
2n-n\leq n$. That is, $|\alpha_i|\leq q^{n/2}$. It follows that
all the eigenvalues $\alpha_i$ must have absolute value $q^{n/2}$.
This completes the proof of Theorem 1.4.
\section{Applications to other situations}
Let $J$ be a subset of $\{1,\cdots,n\}$. For
$\{j_1,\cdots,j_s\}\subseteq J$, we write
$$B_{\{j_1,\cdots,j_s\}}
=\{\sum\limits_{u\in L(\Delta)}a_ux^u\in
B:u_{j_1},\cdots,u_{j_s}>0\}.$$ For $l=0,1,\cdots,n$, we
define$$K_l(f,J)=\bigoplus\limits_{1\leq i_1<\cdots<i_l\leq n}
B_{J\setminus \{i_1,\cdots,i_l\}}e_{i_1}\wedge\cdots\wedge
e_{i_l}.$$ Then $(K_{\bullet}(f,J),\hat{\partial})$ is a
subcomplex of $(K_{\bullet}(f,\emptyset),\hat{\partial})$. The
latter is the complex $(K_{\bullet},\hat{\partial})$ we defined
earlier.
\begin{lemma}The sequence
$$0\rightarrow K_{\bullet}(f,J)
\rightarrow K_{\bullet}(f,J\setminus\{j\})\rightarrow
K_{\bullet}(f^{\{j\}},J\setminus\{j\})\rightarrow0$$is exact,
where $f^{\{j\}}$ is the Witt vector whose $i$-th coordinate is
the sum of monomials of the $i$-th coordinate of $f$ not divided
by $x_j$.
\end{lemma}

We define, for $k=1,2,\cdots$,
$$S_k(f,J) =\sum\limits_{x^{q^k}=x,x_{i_1}\cdots x_{i_r}\neq0} \psi({\text
Tr}_{\mathbb{Q}_p[\mu_{q^k-1}]/\mathbb{Q}_p}(\omega(f)(x)))$$ if
$\{1,\cdots,n\}\setminus J=\{i_1,\cdots,i_r\}$, and $f\in
W_m(\mathbb{F}_q[x_1,\cdots,x_n,(x_{i_1}\cdots x_{i_r})^{-1}])$.
Here the equation $x^{q^k}=x$ is solved in
$(\overline{\mathbb{Q}}_p)^n$. We write
$$L_{f,J}(t)=\exp(\sum\limits_{k=1}^\infty
S_k(f,J)\frac{t^k}{k}).$$ By the above lemma we infer the
following trace formula from the earlier one.
\begin{prop}For
$k=1,2,\cdots$, we have
$$S_k(f,J)=\sum\limits_{l=0}^n (-1)^{l+1}
\text{Tr}((p^nF^{-1})^{ak};H_l(K_{\bullet}(f,J),\hat{\partial})).$$
Equivalently, $$L_{f,J}(t)=\prod\limits_{l=0}^n
\det(1-(p^nF^{-1})^at;H_l(K_{\bullet}(f,J),\hat{\partial}))^{(-1)^l}.$$
\end{prop}

We call $f$ commode with respect to $J$ if $\Delta_{\infty}(f)$ is
commode with respect to $J$. Recall that a convex polyhedron
$\Delta$ in $\mathbb{Q}^n$ that contains the origin is commode
with respect to $J$ if it lies in $(\prod\limits_{i=1,i\not\in
J}^n\mathbb{Q})\times(\prod\limits_{i\in J}\mathbb{Q}_{\geq0})$
and $\dim(\Delta_C)=n-|C|$ for all subset $C$ of $J$, where
$\Delta_C=\{(u_1,\ldots,u_n)\in\Delta:u_j=0\text{ if }j\in C\}$.
By Lemma 11.1 and Proposition 11.2, we infer the following
proposition from Theorem 1.3.
\begin{prop}If $f$ is commode with respect to
$J$ and non-degenerate with respect to $\Delta_{\infty}(f)$, then
$L_{f,J}(t)$ is a polynomial, and its Newton polygon with respect
to $\text{ord}_q$ lies above the Hodge polygon of
$$\sum\limits_{C\subset
J}(-1)^{|C|}P_{\Delta_C}(t^{\frac{D(\Delta)}{D(\Delta_C)}})$$ with
the same endpoints. In particular, $L_{f,J}(t)$ is of degree
$$\sum\limits_{C\subset
J}(-1)^{|C|}(n-|C|)!\text{Vol}(\Delta_C).$$
\end{prop}
From Lemma 10.1 we infer the following one.
\begin{lemma}Let
$\Delta$ be a convex polyhedron in $\mathbb{Q}^n$ of dimension $n$
that contains the origin and is commode with respect to $J$. Let
$(V_{\Delta,J},U_{\Delta,J})$ be the endpoint of the Hodge polygon
of
$$\sum\limits_{C\subset
J}(-1)^{|C|}P_{\Delta_C}(t^{\frac{D(\Delta)}{D(\Delta_C)}})$$
other than $(0,0)$. Then
$$U_{\Delta,J}=\frac
n2V_{\Delta,J}+\sum\limits_{l=1}^{|J|+1}(-1)^l\frac{(n-l)!}{2}(\sum\limits_{C\subset
J,|C|=l-1}S_{\Delta_C}-l\sum\limits_{C\subset
J,|C|=l}\text{Vol}(\Delta_C)).$$ In particular,
$U_{\Delta,J}=\frac n2V_{\Delta,J}$ if the origin is an interior
point of $\Delta_J$.
\end{lemma}
By Lemma 10.3 we infer the following proposition from Proposition
10.2.
\begin{prop}If $f$
is commode with respect to $J$ and non-degenerate with respect to
$\Delta:=\Delta_{\infty}(f)$, and the origin lies in the interior
of $\Delta_J$, then the reciprocal roots of $L_{f,J}(t)$ are of
absolue value $q^{n/2}$.
\end{prop}

\end{document}